\documentclass{amsart}

\usepackage{amssymb}
\usepackage{amsmath}
\input arrow

\renewcommand{\ker}{\textrm{Ker }} \renewcommand{\hom}{\textrm{Hom}}
\DeclareMathOperator{\im}{Im} \DeclareMathOperator{\ext}{Ext} %%Imagen
\DeclareMathOperator{\ch}{Ch}
\DeclareMathOperator{\Mod}{\mbox{-Mod}}
\DeclareMathOperator{\modu}{\mbox{-mod}}
\DeclareMathOperator{\pd}{pd}

\DeclareMathOperator{\cofdim}{cofdim}

\DeclareMathOperator{\id}{id}

\DeclareMathOperator{\Findim}{Findim}

\DeclareMathOperator{\findim}{findim}

\DeclareMathOperator{\fibdim}{fibdim}

\newcommand{\dgF}{{\mathit dg\,}\widetilde{{\mathcal F}}}
\newcommand{\barF}{{\widetilde{{\mathcal F}}}}

\newcommand{\dgC}{{\mathit dg\,}\widetilde{{\mathcal C}}}
\newcommand{\barC}{{\widetilde{{\mathcal C}}}}

\newcommand{\rest}{\upharpoonright}

\newcommand{\fin}{\hspace{\stretch{1}}$\square$}

\newenvironment{demostracion}{\noindent \textbf{Proof.}}

\newtheorem{definicion}{Definition}[section]
\newtheorem{proposicion}[definicion]{Proposition}
\newtheorem{teorema}[definicion]{Theorem}

\newtheorem{lema}[definicion]{Lemma}
\newtheorem{corolario}[definicion]{Corollary}
\newtheorem{remark}[definicion]{Remark}

\title[A MODEL STRUCTURE APPROACH TO THE FINITISTIC DIMENSION CONJECTUREs]{A QUILLEN MODEL STRUCTURE APPROACH\\ TO THE FINITISTIC DIMENSION CONJECTURES}
\author{M. Cort\'es Izurdiaga}

\address{Departamento de \'Algebra y A. Matem\'atico, Universidad de Almer\'{\i}a, Almer\'{\i}a, 04071 } \email{mizurdia@ual.es}

\author{S. Estrada}
\address{
Departamento de Matem\'atica Aplicada, Universidad de Murcia,
Campus del Espinardo, Espinardo (Murcia) 30100, Spain}
\email{sestrada@um.es}
\author{P. A. Guil Asensio}
\address{Departamento de Matem\'aticas, Universidad de Murcia,
Campus del Espinardo, Espinardo (Murcia) 30100, Spain}
\email{paguil@um.es}

 \keywords{projective dimension, model structure, finitistic conjecture.
\\ The first two authors are partially supported by the DGI MTM2008-03339 and by the Junta de Andaluc\'ia FQM-03128.
\\ The last two authors are partially supported by the Fundaci\'on S\'eneca 07552/GERM/2007.}

\begin{document}

\begin{abstract}
  We explore the interlacing between model category structures
  attained to classes of modules of finite $\mathcal{X}$-dimension,
  for certain classes of modules $\mathcal{X}$. As an application we
  give a model structure approach to the Finitistic Dimension
  Conjectures and present a new conceptual framework in which these
  conjectures can be studied.  

\end{abstract}
\maketitle

Let $\Lambda$ be a finite dimensional algebra over a field $k$ (or
more generally, let $\Lambda$ be an artin ring). The big finitistic dimension of
$\Lambda$, $\Findim(\Lambda)$, is defined as the supremum of the
projective dimensions of all modules having finite projective
dimension. And the little finitistic dimension of $\Lambda$,
$\findim(\Lambda)$, is defined in a similar way by restricting to the subclass
of all finitely generated modules of finite projective
dimension. In 1960, Bass stated the so-called Finitistic Dimension
Conjectures: (I) $\Findim(\Lambda) = \findim(\Lambda)$, and (II)
$\findim(\Lambda)$ is finite. The first conjecture was proved to be
false by Huisgen-Zimmermann in \cite{z}, but the second one still
remains open.  It has been proved to be true, for instance, for
finite-dimensional monomial algebras \cite{gkk}, for Artin algebras
with vanishing cube radical \cite{huis}, or Artin algebras with
representation dimension bounded by 3 \cite{it}.

In \cite{a} (see also \cite{HuisSmalo}), Auslander proved that the
finitistic dimension conjectures hold for Artin algebras in which the
class $\mathcal{P}^{<\infty}$ of all finitely generated modules of
finite projective dimension is contravariantly finite (equivalently,
it is a precovering class in the sense of \cite{ej,GT}). In general,
$\mathcal{P}^{<\infty}$ does not need to be contravariantly finite,
even for Artin algebras satisfying the finitistic dimension
conjectures. But, as Angeleri-H\"ugel and Trlifaj have noticed in
\cite{at}, it cogenerates a cotorsion pair $(\mathcal{F},\mathcal{C})$
in which the class $\mathcal F$ is precovering in $R\Mod$. By means of this
idea, the authors are able to extend Auslander's approach to arbitrary
artinian rings and obtain a general criterium for an artinian ring to
satisfy the finitistic dimension conjectures in terms of Tilting
Theory (see \cite{at}). This type of arguments has also been
recently extended to more general homologies induced by arbitrary
hereditary cotorsion pairs (see \cite{amen}).

On the other hand, Hovey has recently shown in \cite{h2} that there exists a
quite strong relation between the construction of hereditary cotorsion
pairs in module categories and the existence of model structures in
the sense of Quillen in the associated categories of unbounded chain
complexes. Moreover, under very general hypotheses, the cohomology
functors defined from these model structures coincide with the absolute
cohomology functors defined from the injective model structure (in the sense
of \cite[Example 3.2]{h2})]. Recall that a model category is a
category with three distinguished classes of morphisms (fibrations,
cofibrations and weak equivalences) satisfying a certain number of
axioms. We refer to \cite{h} for a complete definition and main
properties of model categories. One of the main advantages of these
model categories is that they allow the construction of the derived
category of a ring $R$ as the homotopy category induced by the model
structure. This is the approach followed, for instance, in \cite{eeg},
\cite{egpt}, \cite{g2} and \cite{g1} in order to construct derived
categories of Grothendieck categories in absence of enough projective objects
(in particular, for the category of quasi-coherent sheaves over a
scheme $X$).

The main goal of this paper is to give a new conceptual framework in which
the above results concerning the
finitistic dimension conjectures can be obtained. This is done by developing
Quillen model structures associated to distinguished classes consisting of
modules of finite projective dimension.
In particular, we show that, given any ring $R$, the cotorsion
pair cogenerated by the class of all modules of finite projective
dimension induces a Quillen Model Structure in the category $\ch(R)$
of all unbounded complexes of left $R$-modules,
in which the weak equivalences do coincide with the usual homology
isomorphisms. In particular, by means of this model structure, the absolute
cohomology functors ${\rm Ext}^n(M,N)$ can be recovered in terms of certain
resolutions of $M$ and $N$ attained to the class of modules of finite projective
dimension. This new approach provides a
very general framework in which the different approaches given in
\cite{amen,at,a,t2} find their natural setting. We show that, esentially, they correspond to compute the finitistic dimensions of the considered Artin algebra in a different (and more convenient) homologically equivalent model structure.
We would like to stress that our approach works for any ring $R$, whereas its main
interest appears in the particular case of Artin algebras.

\section{Homology relative to a hereditary cotorsion pair}
\label{sec:preliminaries}

We begin by fixing some notation and terminology. Given a set $X$, we
are going to denote its cardinality by $|X|$; and by $\omega$, the
first infinite ordinal. The {\em cofinality} of an ordinal number $\alpha$
will by denoted by $\textrm{cf}(\alpha)$. I.e., the least
cardinal number which is cofinal in $\alpha$. Recall that an ordinal number
is called {\em regular} when it coincides with its cofinality (and therefore,
it is a cardinal).  The symbol $\upharpoonright$ will mean restriction.

Along this paper, $R$ will denote a ring with identity and all modules
will be left $R$-modules. We will denote by $R\Mod$ the category of all left $R$-modules and by $R\modu$ the subcategory of all modules possessing a projective resolution consisting of finitely generated modules. Morphisms will operate on the right and the
composition of $f:A \rightarrow B$ and $g:B \rightarrow C$ will be
denoted by $fg$. Fixed an infinite regular cardinal $\kappa$, a module
$M$ is said to be {\em $<\kappa$-presented} if it has a free presentation
with less than $\kappa$ generators and relations. If $\mathcal X$ is a
nonempty class of modules, {\em $\mathcal X^{<\kappa}$ will denote the
class of all $<\kappa$-presented modules of $\mathcal X$}.

Let $\mathcal X$ be a nonempty class of modules. We shall consider the
$\ext$-orthogonal classes
\begin{displaymath}
  \mathcal X^\perp = \{Y \in R\Mod: \ext^1(X,Y) = 0 \quad \forall X \in \mathcal{X}\}
\end{displaymath}
and
\begin{displaymath}
  {^\perp}\mathcal{X} = \{Y \in R\Mod: \ext^1(Y,X) = 0 \quad \forall X \in
  \mathcal{X}\}.
\end{displaymath}
Recall that a module {\em $M$ is called $\mathcal X$-filtered} if there
exists, for some regular cardinal $\kappa$, a continuous chain of
submodules of $M$, $\{M_\alpha: \alpha < \kappa\}$ satisfying that the modules $M_0$ and
$\frac{M_{\alpha+1}}{M_\alpha}$ are isomorphic to modules in $\mathcal
X$, for each $\alpha < \kappa$, and $M =
\cup_{\alpha < \kappa}M_\alpha$. An {\em $\mathcal X$-precover} of a module
$M$ is a morphism $\varphi:X \rightarrow M$ such that $X \in \mathcal
X$ and $\hom(X',\varphi)$ is an epimorphism for every $X' \in \mathcal
X$. An $\mathcal X$-precover $\varphi:X \rightarrow M$ of $M$ is
called {\em special} if it is an epimorphism and $\ker \varphi \in \mathcal
X^\perp$ (see e.g., \cite{GT}). $\mathcal X$-preenvelopes and special $\mathcal
X$-preenvelopes are defined dually.

A {\em cotorsion pair} in $R\Mod$ is a pair of classes of modules,
$(\mathcal F, \mathcal C)$, such that $\mathcal F = {}^{\perp}\mathcal
C$ and $\mathcal C = \mathcal F^\perp$. The cotorsion pair is said to
be {\em cogenerated by a class of modules $\mathcal{X}$} if
$\mathcal{X}^\perp = \mathcal C$. When this class $\mathcal{X}$ is a set, it is
known that every module has a special $\mathcal F$-precover and a
special $\mathcal C$-preenvelope (see e.g. \cite[Theorem
3.2.1]{GT}). In this case the cotorsion pair is called complete.

Let $(\mathcal F,\mathcal C)$ be a cotorsion pair cogenerated by a
set. Then there exists an infinite regular cardinal
$\kappa$ such that $(\mathcal F,\mathcal C) $ is cogenerated by
$\mathcal F^{<\kappa}$. Moreover, by Kaplansky's Theorem for cotorsion
pairs (see \cite[Theorem 10]{st} or \cite[Theorem 4.2.11]{GT}),
$\mathcal F$ consists of all $\mathcal F^{<\kappa}$-filtered modules.

A cotorsion pair $(\mathcal F, \mathcal C)$ is called {\em hereditary} if
the class $\mathcal{F}$ is {\em resolving}. I. e., it is closed under
kernels of surjections and contains all projective modules.
%This is
%equivalent to say that $\mathcal C$ is coresolving (i. e., it is
%closed under cokernels of injections and contains all injective
%modules).

We now recall some well-known facts of the category $\textrm{Ch}(R)$
of unbounded chain complexes of modules. A complex of $R$-modules,
\begin{displaymath}
  \harrowlength=30pt \cdots \mapright^{d_{n+2}} X_{n+1}
  \mapright^{d_{n+1}} X_n \mapright^{d_n} X_{n-1}
  \mapright^{d_{n-1}} \cdots,
\end{displaymath}
will be denoted by $(X,d)$, or simply by $X$. And we will denote
by $Z_nX = \ker d_n$, $K_nX = \textrm{Coker } d_n$, $B_nX
= \im d_{n+1}$ and $H_nX = \frac{Z_nX}{B_nX}$, for every integer
$n$. Given other complex $Y$, $Hom(X,Y)$ will denote the complex
defined by
\begin{displaymath}
  Hom(X,Y)_n = \prod_{k \in \mathbb{Z}}{\rm Hom}_R(X_k,Y_{k+n})
\end{displaymath}
and $\big((f)d^H_n\big)_k  = f_kd_{k+n}^Y-(-1)^nd_k^Xf_{k-1}$
for any $n \in \mathbb{Z}$. The class of all exact complexes will be
denoted by $\mathcal E$.

Let us fix a cotorsion pair $(\mathcal{F}, \mathcal{C})$ in
$R\Mod$. We will consider the following subclasses of $\textrm{Ch}(R)$
(see \cite[Definition 3.3]{g1}):

\begin{enumerate}
\item The class of {\em $\mathcal F$-complexes}, $\barF = \{X \in
  \mathcal{E}:Z_nX \in \mathcal{F}, \ \forall n \in \mathbb{Z}\}$.

\item The class of {\em $\mathcal C$-complexes}, $\barC = \{X \in
  \mathcal{E}:Z_nX \in \mathcal{C}, \ \forall n \in \mathbb{Z}\}$.

\item The class of {\em dg-$\mathcal F$ complexes},
  $$\dgF = \{X \in \ch(R): X_n \in
  \mathcal{F} \ \forall n \in \mathbb{Z} \textrm{ and } Hom(X,C)
  \textrm{ is exact}\ \forall C \in \tilde{\mathcal{C}}\}.$$

\item The class of {\em dg-$\mathcal C$ complexes},
  $$\dgC = \{X \in \ch(R): X_n \in
  \mathcal{C}\ \forall n \in \mathbb{Z} \textrm{ and } Hom(F,X)
  \textrm{ is exact} \ \forall F \in \tilde{\mathcal{F}}\}.$$
\end{enumerate}

Our next theorem shows that any hereditary cotorsion pair $(\mathcal
F, \mathcal C)$ in $R\Mod$ cogenerated by a set gives rise to a model
structure in $\textrm{Ch}(R)$ in which the weak equivalences are the
homology isomorphisms. This is essentially due to Hovey \cite[Theorem
2.2]{h2}.

\begin{teorema}\label{t:ModelStructure}
  Let $(\mathcal{F},\mathcal{C})$ be a hereditary cotorsion pair in
  $R\Mod$ cogenerated by a set. Then there exists a cofibrantly
  generated model structure on $\ch(R)$ such that:
  \begin{enumerate}
  \item the weak equivalences are the homology isomorphisms;

  \item the cofibrations (resp. trivial cofibrations) are the
    monomorphisms whose cokernels are in
    $\textrm{dg\,}\tilde{\mathcal{F}}$ (resp. $\tilde{\mathcal{F}}$); and

  \item the fibrations (resp. trivial fibrations) are the
    epimorphisms whose kernels are in
    $\textrm{dg\,}\tilde{\mathcal{C}}$ (resp. $\tilde{\mathcal{C}}$).
  \end{enumerate}
\end{teorema}

\begin{demostracion}
  We are going to apply \cite[Theorem 2.2]{h2}. So we must check the
  following conditions:
  \begin{enumerate}
  \item The pairs $(\barF,\dgC)$ and $(\dgF,\barC)$ are cotorsion
    pairs.
  \item $\dgF\cap \mathcal E=\barF$ and $\dgC\cap \mathcal E=\barC$.
  \item The cotorsion pairs $(\barF,\dgC)$ and $(\dgF,\barC)$ are
    complete.
  \end{enumerate}
  By \cite[Corollary 3.8]{g1} we have induced cotorsion pairs
  $(\barF,\dgC)$ and $(\dgF,\barC)$. Now (2) follows from
  \cite[Corollary 3.12]{g1}.  By \cite[Proposition 3.8]{g2} and \cite[Corollary
  6.6]{h2} the pair $(\dgF,\barC)$ is complete. Finally to see that the pair $(\barF,\dgC)$ is complete we
  appeal to the same proof of \cite[Lemma 4.2]{egpt} (notice that the
  arguments of the proof of \cite[Lemma 4.2]{egpt} are categorical, so
  they can be adapted to our setting).
  \fin \end{demostracion}
\begin{remark}
  Let us note that we cannot use \cite[Theorem 4.12]{g2} to get the desired model
  structure since we are not assuming that the class $\mathcal F$ is
  closed under direct limits. We cannot assume this closure under direct limits
  since this condition is rarely satisfied in the applications we are interested in.
\end{remark}
We are going to denote by $\ch(R)_{\mathcal{M}_\mathcal{F}^\mathcal{C}}$ the
category $\ch(R)$ endowed with the above model structure induced by the
class $\mathcal{F}$. Let us recall that, if ${\rm Ho}(\ch(R))$ is the homotopy category associated to $\ch(R)_{\mathcal{M}_\mathcal{F}^\mathcal{C}}$, we can define
$\ext^n_R(M,N)$ as $${\rm Hom}_{{\rm Ho}(\ch(R))}(S(M),S^n(N))={\rm
  Hom}_{\ch(R)}(Q_{S(M)},P_{S^n(N)})/\sim_h,$$(see \cite[Theorem 1.2.10]{h}) where $Q_{S(M)}$ is a
cofibrant replacement of $S(M)$ (i.e., the complex with $M$ in the 0'th position
and 0 elsewhere), $P_{S^n(N)}$ is a fibrant replacement of $S^n(N)$
(the complex with $N$ in the $n$'th position and 0 elsewhere) and $\sim_h$
denotes the chain homotopy (see \cite{h}).

Given a module $M$, the standard way of constructing a cofibrant replacement $Q_M$
of $S(M)$ (that is, a trivial fibration $Q_M\to S(M)$, where $Q_M$ is
cofibrant) from a hereditary cotorsion pair $(\mathcal F,\mathcal C)$ as in Theorem \ref{t:ModelStructure}, is the following: we choose a special $\mathcal F$-precover $d_0: F_0\to M\to 0$
 of $M$. Then, a
special $\mathcal F$-precover $d_1:F_1\to \ker d_0\to 0$ of $\ker d_0$ and so on.
Proceeding
in this way, we get an exact complex $\cdots
\stackrel{{d_2}}{\mapright} F_1
\stackrel{{d_1}}{\mapright} F_0
\stackrel{{d_0}}{\mapright} M \to 0,$ in which $F_i \in \mathcal
F$ and $\ker d_i \in \mathcal F^\perp=\mathcal C$ for every $i \in \mathbb{N}$. Then,
if we denote by $F_\bullet$ (or by $(F_\bullet, d)$) the
corresponding deleted complex (which is unique up to chain homotopy
equivalence), we get an epimorphism in $\ch(R)$, $\bar{d}:F_\bullet\to
S(M)$, with $\ker(\bar{d})\in \tilde{\mathcal{C}}$ (and $F_\bullet\in
dg\,\tilde{\mathcal{F}}$) and therefore, a cofibrant replacement of
$S(M)$. Dually, we can get a fibrant replacement of $S(N)$,
$C_{\bullet}$ (unique up to chain homotopy equivalence), from the fact
that every module admits a special $\mathcal C$-preenvelope. Notice
that both fibrant and cofibrant replacements are unique in
${\rm Ho}(\ch(R))$ (because they provide unique-up-to-homotopy
resolutions and coresolutions, respectively).
Then, as $${\rm
  Hom}_{\ch(R)}(Q_{S(M)},P_{S^n(N)})/\sim_h\cong
H_{-n}Hom(F_\bullet,C_\bullet),$$  we may identify
the Ext groups $\ext^n_R(M,N)$ with the homology groups of the complex
$Hom(F_\bullet,C_\bullet)$. Using this identification, we can easily
deduce the following result.

\begin{proposicion}\label{p:CalculatingExt}
  Let $(\mathcal{F},\mathcal{C})$ be a hereditary cotorsion pair
  cogenerated by a set, and $M,N$, two modules.
  \begin{enumerate}
  \item If $M \in \mathcal{F}$ and $C_\bullet$ is a fibrant
    replacement of $S(N)$, then:
    \begin{enumerate}
    \item $\ext^n(M,N) = H_{-n}Hom(M,C_\bullet)$.

    \item $\ext^n(M,K_{-k}C_\bullet) = \ext^{n+k+1}(M,N)$ for every $k
      \leq n$.
    \end{enumerate}

  \item If $N \in \mathcal{C}$ and $F_\bullet$ is a cofibrant
    replacement of $S(N)$, then:
    \begin{enumerate}
    \item $\ext^n(M,N) = H_{-n}Hom(F_\bullet,N)$.

    \item $\ext^n(Z_{k}F_\bullet,N) = \ext^{n+k+1}(M,N)$ for every $k
      \geq n$, .
    \end{enumerate}
  \item If $M \in \mathcal F$ and $N \in \mathcal C$, then
    $\ext^n(M,N) = 0$ for every $n \geq 1$.
  \end{enumerate}
\end{proposicion}

The above results suggest the definition of some homological invariants with respect
to the classes $\mathcal F$
and $\mathcal C$. Let $X$ be a complex and $a,b \in \mathbb{Z}$ with
$a < b$. We will say that the amplitude of $X$ belongs to $[a,b]$ if $X_n
= 0$ for every $n < a$ and every $n > b$.

\begin{definicion}
  Let $(\mathcal{F},\mathcal{C})$ be a hereditary cotorsion pair
  cogenerated by a set, and $M,N$, two modules. Then:
  \begin{enumerate}
  \item Given an $n\in \mathbb{N}$, we will say that the
    $\mathcal{F}$-cofibration dimension of $M$ is bounded by $n$ if
    there exists a cofibrant replacement $F_\bullet$ of $S(M)$ in
    $\ch(R)_{\mathcal{M}_\mathcal{F}^\mathcal{C}}$ with amplitude in
    $[0,n]$.  We will define the $\mathcal{F}$-cofibrant dimension of
    $M$, denoted by $\cofdim_{\mathcal{F}}(M)$, as the minimal natural
    number $n$ satisfying this property and $\infty$ otherwise.
  \item Given an $n\in \mathbb{N}$, we will say that the projective
    dimension of $M$ relative to $\mathcal{C}$ is bounded by $n$ if
    $\ext^m(M,C) = 0$ for every $m > n$ and $C \in \mathcal C$. And
    we will define the projective dimension of $M$ relative to
    $\mathcal{C}$, $\pd_\mathcal{C}(M)$, as the minimal $n\in \mathbb{N}$
    satisfying this property and $\infty$ otherwise.
  \end{enumerate}

  Analogously we may define the injective dimension of $M$ relative to
  $\mathcal F$, $\id_\mathcal{F}(M)$, and the $\mathcal{C}$-fibration
  dimension of $M$, $\fibdim_\mathcal{C}(M)$.
\end{definicion}

Let $\mathcal A$ be any Grothendieck category. If we consider the
injective model structure in $\ch(\mathcal A)$ introduced by Joyal
in \cite{joyal} and Beke in \cite{Beke} (see also \cite[Theorem 2.3.13]{h}),
it is easy to check that the
fibrant dimension of an object $M$ of $\mathcal A$ is precisely the
usual injective dimension of $M$ (note that, in this case, a fibrant replacement
of $M$ is an injective coresolution of $M$). Analogously, Hovey has proved in
\cite[Chapter 4]{h} that there exists a projective monoidal
model structure in $\ch(R)$. And with respect to this model structure,
the cofibrant dimension of a module $M$ is precisely the usual
projective dimension of $M$. We can obtain these two model structures
in $\ch(R)$ from Theorem \ref{t:ModelStructure}, by considering the
hereditary cotorsion pairs $(\mathcal Proj, R\Mod)$ and $(R\Mod,
\mathcal Inj)$, which are obviously cogenerated by sets. We are going to extend in
Proposition \ref{p:CalculatingDimensions} the usual properties of
the classical injective and projective dimension to the dimensions
induced from any hereditary and complete cotorsion pair. In order to prove it,
we will need the following generalized version of Schanuel's lemma (see
\cite[Lemma 8.6.3]{ej} for a proof):
\begin{lema}\label{l:shortexactsequence}
  Let $\mathcal F$ be a class of modules, $M$ a module, $F,G \in
  \mathcal F$, $\varphi:F \rightarrow M$ a surjective $\mathcal
  F$-precover and $f:G \rightarrow M$, an epimorphism. Then there
  exists a short exact sequence $0\to \ker f\to \ker \varphi \oplus  G\to F\to 0 $.
\end{lema}
\begin{proposicion}\label{p:CalculatingDimensions}
  Let $(\mathcal F, \mathcal C)$ be an hereditary cotorsion pair
  cogenerated by a set, $M$, a module and $n$, a natural number. The
  following assertions are equivalent:
    \begin{enumerate}
    \item $\textrm{cofdim}_\mathcal{F}(M) \leq n$.
    \item $\pd_\mathcal{C}(M) \leq n$.
    \item Each cofibrant replacement $Q$ of $S(M)$ satisfies that $Z_{n-1}Q
      \in \mathcal F$.
    \item There exists a short exact sequence
      \begin{displaymath}
        \harrowlength=20pt 0 \mapright G_n \mapright^{d_n} \cdots \mapright G_0
        \mapright^{d_0} M \mapright 0
      \end{displaymath}
      with $G_i \in \mathcal F$ for every $i$.
    \end{enumerate}
\end{proposicion}
\begin{demostracion}
    (1) $\Leftrightarrow$ (2) follows from Proposition
    \ref{p:CalculatingExt}. (1) $\Leftrightarrow$ (3) is clear since we can
    compute Ext functors by using any cofibrant replacement of
    $S(M)$. (3) $\Rightarrow$ (4) is trivial.
In order to finish the proof we need to show that
    (4) $\Rightarrow$ (3). Let $Q$ be a cofibrant replacement of
    $S(M)$. We can assume that $Q$ is bounded below, and thus it is of the
    form $(F_\bullet,\varphi)$. We will induct on $n\in\mathbb{N}$. If $n
    =0$, then $M \in \mathcal F$ and the result is clear.  Suppose now
    that the result is true for any module $M$ with
    $\cofdim_\mathcal{F}(M)\leq n$ and let us prove it for
    $\cofdim_\mathcal{F}(M)\leq n+1$. Then there exists an exact
    sequence with $n+2$ terms, $ 0 \to G_{n+1}
    \stackrel{d_{n+1}}{\mapright}G_{n} \to \ldots \to G_1
    \stackrel{d_1}{\mapright} G_0
    \stackrel{d_0}{\mapright} M \to 0$ with $G_i \in \mathcal F$
    for every $i$. Let us denote by $\psi_1: F_1 \oplus G_0
    \rightarrow \ker \varphi_0 \oplus G_0$ the direct sum $\varphi_1
    \oplus 1_{G_0}$, and let us note that it is a special $\mathcal
    F$-precover. Then, if we construct the pullback of the short exact
    sequence obtained in Lemma \ref{l:shortexactsequence} for the
    epimorphism $d_0$ and the $\mathcal F$-precover $\varphi_0$, we
    get the following commutative diagram with exact rows and columns:
  \begin{displaymath}
    \commdiag{ & & 0 & & 0 & & & & \cr
      & & \downarrow & & \downarrow & & & & \cr
      & & \ker \varphi_1 & \mapright & \ker \varphi_1 & & & & \cr
      & & \downarrow\lft{\overline{\iota}_1} & & \downarrow\rt{\iota_1} & & & & \cr
      0 & \mapright & P & \mapright^{\overline{\alpha}} &
      F_1\oplus G_0 &
      \mapright^{\psi_1 \beta} & F_0 & \mapright & 0\cr
      & & \downarrow\lft{\overline{\psi_1}} & &
      \downarrow\rt{\psi_1} & & \downarrow & & \cr
      0 & \mapright & \ker d_0 & \mapright_\alpha & \ker
      \varphi_0 \oplus G_0 &
      \mapright_\beta & F_0 & \mapright & 0\cr
      & & \downarrow & & \downarrow & & & & \cr
      & & 0 & & 0 & & & &}.
  \end{displaymath}
  Now note that, since $\mathcal F$ is resolving, $P \in \mathcal
  F$. And, as the bottom square is a pullback, $\overline{\psi}_1$ is
  a special $\mathcal F$-precover. Therefore, $\ker d_0$ has a
  cofibrant replacement $Q'$ given by the deleted complex of $ \cdots \to
  F_3 \stackrel{{\varphi_3}}{\mapright} F_2
  \stackrel{{\varphi_2\overline{\iota}_1}}{\mapright} P
  \stackrel{{\overline{\psi}_1}}{\mapright} \ker d_0 \to 0.$
  Finally, as $\cofdim_\mathcal{F}(\ker d_0)\leq n$, we get from the induction hyphotesis that $K_{n-1}Q' \in
  \mathcal F$. But $K_{n}Q = K_{n-1}Q'$, and the result is proved.\fin
  \end{demostracion}

As a consequence of this result we get:

\begin{corolario}
  Let $(\mathcal F, \mathcal C)$ be a hereditary cotorsion pair
  cogenerated by a set and $M$, a module. Then
  $\cofdim_{\mathcal{F}}(M) = \pd_\mathcal{C}(M)$ and
  $\fibdim_{\mathcal{C}}(M) = \id_\mathcal{F}(M)$.
\end{corolario}

\section[Model structures from finite relative dimensional
modules]{Model structures from finite relative dimensional modules}
\label{sec:model-struct-induc}

Throughout this section, we will fix a hereditary cotorsion pair
$(\mathcal{F},\mathcal{C})$ cogenerated by a set. For any $n \in
\mathbb{N}$, we will denote by $\mathcal P_n$ (resp. $\mathcal
P_n^{<\infty}$) the class of all modules (resp. all modules in $R\modu$)
having projective dimension relative to $\mathcal C$ at most $n$. And by
$\mathcal P$ (resp.  $\mathcal P^{<\infty}$) the class $\bigcup_{n \in
  \mathbb{N}} \mathcal P_n$ (resp. $\bigcup_{n\in \mathbb{N}} \mathcal
P_n^{<\infty}$). Note that, by Proposition \ref{p:CalculatingDimensions},
$\mathcal P_n$ coincides with the class of all modules $M$ having an exact sequence
$0\to F_j\to F_{j-1}\to \ldots\to F_0\to M\to 0$ with $F_i\in \mathcal F$, for
any $1\leq i\leq j$ and any $j\leq n$.
Finally, $(\mathcal A_n, \mathcal B_n)$ (resp
$(\mathcal A_n^{<\infty},\mathcal B_n^{<\infty})$) will denote the
cotorsion pair cogenerated by $\mathcal P_n$ (resp.  $\mathcal
P_n^{<\infty}$), and $(\mathcal A, \mathcal B)$ (resp. $(\mathcal
A^{<\infty}, \mathcal B^{<\infty})$), the cotorsion pair cogenerated
by $\mathcal P$ (resp.  $\mathcal P^{<\infty}$).

The following lemma generalizes \cite[Proposition 3]{a} (just take $\mathcal C=R\Mod$).
\begin{lema}\label{closure}
The class $\mathcal P_n=\{L: \pd_{\mathcal C}(L)\leq n\}$ is closed under $\mathcal P_n$-filtrations for any $n\geq0$.
\end{lema}
\begin{demostracion}
Let $\mathcal D$ be the class of all $k$-th cosyzygies of objects of $\mathcal C$ ($k\geq n$). Then a module
$L\in {}^{\perp} \mathcal D$ if and only if ${\rm Ext}^{m}(L,C)=0$, for each $C\in \mathcal C$ and each $m>n$.
And this happens if and only if $L\in \mathcal P_n$. By \cite[Theorem 1.2]{e} it follows that $\mathcal P_n$ is closed under $\mathcal P_n$-filtrations.
\fin
\end{demostracion}

%\begin{remark}
%Although $\mathcal A_n=\mathcal P_n$ we prefer to keep the notation of $\mathcal %A_n$ for the class of modules with finite projective dimension relative to %$\mathcal C$, because $\mathcal P_n$ is commonly used to denote the class of %modules with finite projective dimension $\leq n$.
%\end{remark}

The main goal of this section will be to show that the above cotorsion pairs
induce model structures in $\ch(R)$. By Theorem \ref{t:ModelStructure}, we only
need to prove that they are hereditary and cogenerated by a set. Let
us note that the cotorsion pairs $(\mathcal A^{<\infty},
\mathcal B^{<\infty})$ and $(\mathcal A^{<\infty}_n,\mathcal
B^{<\infty}_n)$ (for each $n\in \mathbb{N}$) are cogenerated by a set by definition. The
next result shows that the same holds for $(\mathcal A, \mathcal B)$ and $(\mathcal A_n,\mathcal
B_n)$ ($n\in \mathbb{N}$).

\begin{teorema}\label{t:cogenerated}
  Let $(\mathcal F, \mathcal C)$ be a hereditary cotorsion pair
  cogenerated by a set.  Then $(\mathcal A, \mathcal B)$ and $(\mathcal A_n,\mathcal B_n)$ are also
cogenerated by a set, for each $n\in \mathbb{N}$.
\end{teorema}

\begin{demostracion}
It suffices to show that $(\mathcal A_n, \mathcal B_n)$ is
cogenerated by a set for each $n\in \mathbb{N}$, since $\mathcal A=\cup_{n\in \mathbb{N}}\mathcal A_n$. Let us fix an $n\in \mathbb{N}$. Since $(\mathcal F,\mathcal C)$ is cogenerated by a set,
 there exists by \cite[Lema 4.2.10]{GT} an infinite regular cardinal $\kappa$ such that each module in $\mathcal F$ is
$\mathcal F^{<\kappa}$-filtered. We are going to prove that each module in $\mathcal P_n$ is $\mathcal P_n^{<\kappa}$-filtered. Then,
if $\mathcal S$ is a set of representatives of the isomorphism classes of modules in $\mathcal P_n^{<\kappa}$, we get from \cite[Theorem 1.2]{e} that the cotorsion pair $(\mathcal A_n,\mathcal B_n)$ is cogenerated by $\mathcal S$.

Let us fix an $A \in \mathcal P_n$ and an
exact sequence,
\begin{equation}\label{eq:1}
    0 \to F_n \stackrel{{d_n}}{\mapright} F_{n-1} \stackrel{{d_{n-1}}}{\to} \ldots
    \stackrel{{d_1}}{\mapright} F_0 \stackrel{{d_0}}{\mapright} A \to 0
\end{equation} with $F_i\in \mathcal F$, $1\leq i\leq n$.
We know that each $F_i$ is $\mathcal F^{<\kappa}$-filtered. Let us denote by $\mathcal F_i$ a family of submodules
of $F_i$ given by the Hill Lemma (see \cite[Theorem 4.2.6]{GT}). We will follow the notation used in \cite[Theorem 4.2.6]{GT}. In particular, we will refer to the properties satisfied by this family as properties (H1), (H2),
  (H3) and (H4).

 Let $\{x_\alpha:\alpha < \mu\}$ be a
  generating set of $A$ with $x_0 = 0$. For each $\alpha < \mu$, we are going to
  construct an exact sequence $S_\alpha$,
  \begin{displaymath}
    0 \to F_n^{\alpha} \stackrel{{d_n^\alpha}}{\mapright}
    F_{n-1}^\alpha
    \stackrel{{d_{n-1}^\alpha}}{\mapright} \ldots
    \stackrel{{d_1^\alpha}}{\mapright}
    F_0^{\alpha} \stackrel{{d_0^\alpha}}{\mapright} A_\alpha \to 0
  \end{displaymath}
  such that $x_\alpha \in A_\alpha$, $d_i^\alpha$ is the restriction of $d_i$,
  $\{F_i^\alpha:\alpha < \mu\}$ is an $\mathcal F^{<\kappa}$-filtration
  of $\bigcup_{\alpha < \mu}F^\alpha_i$ for
  every $i \in \{0, \ldots, n\}$, and $\{A_\alpha:\alpha < \mu\}$ is an $\mathcal
  A_n^{<\kappa}$-filtration of $A$.

  We will make this construction by transfinite induction on $\alpha<\mu$. For $\alpha =
  0$, let us fix $A_0 = 0$. Assume now that $0\leq\alpha<\mu$ and that
  $S_\alpha$ has been already constructed. We are going to find, for
  each $i \in \{0, \ldots, n\}$, two chains of submodules of $F_i$,
  $\{X^{\alpha+1}_{im}: m \in \mathbb{N}\}$ and $\{Y^{\alpha+1}_{im}:m \in
  \mathbb{N}\}$, satisfying the following properties:
  \begin{enumerate}
  \item[a)] $A_\alpha \cup \{x_{\alpha+1}\} \leq
    (X_{00}^{\alpha+1})d_0$ and $F_i^\alpha \leq X_{i0}^{\alpha+1}$
    for each $i \in \{0, \ldots, n\}$.

  \item[b)] $X^{\alpha+1}_{im} \leq Y^{\alpha+1}_{im} \leq
    X^{\alpha+1}_{i,m+1}$ for each $m \in \mathbb{N}$ and each $i \in \{0,
    \ldots, n\}$.

  \item[c)] $\ker \left(d_i
      \upharpoonright X^{\alpha+1}_{im}\right) \leq
    \left(X^{\alpha+1}_{i+1m}\right)d_{i+1}$ for each $i \in \{0,
    \ldots, n-1\}$ and each $m \in \mathbb{N}$.

  \item[d)] $\left(Y^{\alpha+1}_{im}\right)d_i
    \leq Y^{\alpha+1}_{i-1,m}$ for each $i \in \{0, \ldots, n-1\}$
    and each $m \in \mathbb{N}$.

  \item[e)] The quotients $\frac{X_{i0}}{F_i^\alpha}$,
    $\frac{Y_{im}^{\alpha+1}}{X_{im}^{\alpha + 1}}$ and
    $\frac{X_{i,m+1}^{\alpha+1}}{Y_{im}^{\alpha+1}}$ belong to
    $\mathcal F^{<\kappa}$ for each $i \in \{0, \ldots, n\}$ and each $m \in
      \mathbb{N}$.
  \end{enumerate}

  We will proceed by induction on $m\in \mathbb{N}$. If $m = 0$, take $y_{\alpha+1} \in F_0^\alpha$
  with $(y_{\alpha+1})d_0 = x_{\alpha+1}$. Since $F_0^\alpha \in
  \mathcal F_0$ there exists, by (H4), an
  $X^{\alpha+1}_{00} \in \mathcal F_0$ with
  $F_0^\alpha \cup \{y_{\alpha+1}\} \leq
  X^{\alpha+1}_{00}$ and
  $\left|\frac{X_{00}^{\alpha+1}}{F_0^\alpha}\right| < \kappa$. Note
  that, by (H3), this module is $\mathcal F$-filtered and therefore, it belongs to $\mathcal F$ by hypothesis.

  In order to construct $X_{10}^{\alpha+1}$ note that, as
  $\left|\frac{X_{00}^{\alpha+1}}{F_0^\alpha}\right| < \kappa$, we deduce from the Ker-Coker
Lemma that $\left|\frac{\ker d_0 \rest
      X_{00}^{\alpha+1}}{\ker d_0 \rest F_0^\alpha}\right| <
  \kappa$. Thus, there exists a submodule $K \leq F_1$ with cardinality strictly smaller
  than $\kappa$ such that $\ker d_0 \rest X_{00}^{\alpha+1} \leq
  (F_1^\alpha \cup K)d_1$. Using again (H4), we may find an
  $X^{\alpha+1}_{10} \in \mathcal F_1$ with $F_1^\alpha \cup K \leq
  X^{\alpha+1}_{10}$ and
  $\left|\frac{X^{\alpha+1}_{10}}{F_1^\alpha}\right| <
  \kappa$. Moreover, this module is $\mathcal
  F$-filtered by (H3) and therefore, it belongs to $\mathcal F$. Repeating a
  similar argument we may find $X^{\alpha+1}_{20}, \ldots,
  X^{\alpha+1}_{n0}$ satisfying the desired properties.

  Now we proceed to construct the modules $Y's$. Let us take $Y^{\alpha+1}_{n0}
  = X^{\alpha+1}_{n0}$. Using that $\frac{Y^{\alpha+1}_{n0}}{F_n^\alpha}$
  has cardinality strictly smaller than $\kappa$, it is easy to check that
  $\frac{(Y_{n0}^{\alpha+1})d_n+X_{n-1,0}^{\alpha+1}}{X_{n-1,0}^{\alpha+1}}$
  also has cardinality strictly smaller than $\kappa$. So, by (H4), there exists
  an $Y_{n-1,0}^{\alpha+1} \in \mathcal F_{n-1}$ with
  $(Y_{n0}^{\alpha+1})d_n+X_{n-1,0}^{\alpha+1} \leq
  Y_{n-1,0}^{\alpha+1}$ and
  $\left|\frac{Y_{n-1,0}^{\alpha+1}}{X_{n-1,0}^{\alpha+1}}\right| <
  \kappa$. Note that, in addition, this quotient belongs to $\mathcal
  F$ since it is $\mathcal F$-filtered by (H3). Repeating this procedure, we get
  $Y^{\alpha+1}_{n-2,0}, \ldots, Y^{\alpha+1}_{00}$. This finishes the case $m=0$.

  Suppose now that we have constructed the chains for some $m \in
  \mathbb{N}$ and take $X^{\alpha+1}_{0,m+1} =
  Y^{\alpha+1}_{0m}$. Since
  $\left|\frac{X^{\alpha+1}_{0,m+1}}{F_0^\alpha}\right| < \kappa$, then
  $\left|\frac{\ker d_0 \rest
      X_{0,m+1}^{\alpha+1}+(Y_{1,m}^{\alpha+1})d_1}{\left(Y_{1,m}^{\alpha+1}\right)
      d_1}\right|
  < \kappa$ and, by (H4), there exists $X_{1,m+1}^{\alpha+1} \in
  \mathcal F_1$ with
  $\ker d_0 \rest X_{0,m+1}^{\alpha+1}+Y_{1,m}^{\alpha+1} \leq
  X_{1,m+1}^{\alpha+1}$ and
  $\left|\frac{X_{1,m+1}}{Y_{1,m}^{\alpha+1}}\right| < \kappa$. Note
  that, in addition, this module belongs to $\mathcal F$ since it is
  $\mathcal F$-filtered by (H3). Repeating this process we get
  $X_{2,m+1}^{\alpha+1}, \ldots, X_{n,m+1}^{\alpha+1}$.

  In order to construct the $Y's$, set $Y_{n,m+1}^{\alpha+1} =
  X_{n,m+1}^{\alpha+1}$. Using that
  $\left|\frac{Y^{\alpha+1}_{0,m+1}}{F_0^\alpha}\right| < \kappa$ we
  conclude easily that
  $\left|\frac{(Y_{n,m+1}^{\alpha+1})d_n+X_{n-1,m+1}^{\alpha+1}}{X_{n-1,m+1}^{\alpha+1}}\right|
  < \kappa$ and, applying (H4) again, there exists an
  $Y_{n-1,m+1}^{\alpha+1} \in \mathcal F_{n-1}$ with
  $\left(Y_{n,m+1}^{\alpha+1}\right)d_n+X_{n-1,m+1}^{\alpha+1} \leq
  Y_{n-1,m+1}^{\alpha+1}$ and
  $\left|\frac{Y_{n-1,m+1}^{\alpha+1}}{X_{n-1,m+1}^{\alpha+1}}\right|
  < \kappa$. Note that, in addition, this module belongs to $\mathcal
  F$ as it is $\mathcal F$-filtered. Proceeding in this way we construct
  $Y_{n-2,m+1}^{\alpha+1}, \ldots, Y_{0,m+1}^{\alpha+1}$.

  Set $F_i^{\alpha+1} = \bigcup_{m \in
    \mathbb{N}}X_{i,m}^{\alpha+1} = \bigcup_{m \in
    \mathbb{N}}Y_{i,m}^{\alpha+1}$ for each $i \in \{0, \ldots, n\}$
  and note that the sequence $S_{\alpha+1}$ is exact by c) and
  d). Since $S_\alpha$ is a subcomplex of $S_{\alpha+1}$, the
  corresponding quotient complex is an exact sequence with each term,
  except the first one, in $\mathcal F^{<\kappa}$ by d). This means
  that $\frac{A_{\alpha+1}}{A_\alpha} \in \mathcal P_n^{<\kappa}$. This finishes the case $\alpha+1$.

  Finally, suppose that $\alpha$ is a limit
  ordinal and define $S_\alpha$ taking $F_i^\alpha =
  \bigcup_{\gamma<\alpha}F_i^\alpha$ and $A_\alpha = \bigcup_{\gamma <
    \alpha}A_\gamma$. It is easy to see that $S_\alpha$ satisfy the
  desired properties.\fin

\end{demostracion}

\begin{corolario}
Let $(\mathcal F, \mathcal C)$ be a hereditary cotorsion pair
  cogenerated by a set.  Then $\mathcal A_n=\mathcal P_n$ for each $n\in \mathbb{N}$.
\end{corolario}

\begin{demostracion}
As the cotorsion pair $(\mathcal A_n,\mathcal B_n)$ is cogenerated by a set $\mathcal S\subseteq \mathcal P_n$, the result follows from Lemma \ref{closure} and \cite[Corollary 3.2.4]{GT} (notice that $\mathcal P_n$ is closed under direct summands).\fin
\end{demostracion}

Let us now check that the above cotorsion pairs are hereditary. In order to prove it, we
will need to use the following two propositions. The proof of the first one is straightforward.

\begin{proposicion}\label{p:hereditary}
  Let $(\mathcal{F},\mathcal{C})$ be a cotorsion pair. The following
  assertions are equivalent:
  \begin{enumerate}
  \item The cotorsion pair is hereditary.

  \item Every projective presentation
    $\varphi:P \rightarrow F$ of any element $F \in \mathcal{F}$ verifies that $\ker \varphi \in
    \mathcal{F}$.

  \item For every $F \in \mathcal{F}$, there exists a projective
    presentation $\varphi:P \rightarrow F$ with $\ker \varphi \in
    \mathcal{F}$.
  \end{enumerate}
\end{proposicion}

\begin{proposicion}\label{p:filtered}
  Let $\mathcal{X}$ be a class of modules such that for each $X \in
  \mathcal X$, there exits a free presentation $\varphi:R^{(I)}
  \rightarrow X$ with $\ker \varphi \in \mathcal X$. If $A$ is a
  direct summand of an $\mathcal{X}$-filtered module, then there exists a
  free presentation $\psi:R^{(J)} \rightarrow A$ in which $\ker \psi$
  is also a direct summand of an $\mathcal{X}$-filtered module.
\end{proposicion}

\begin{demostracion}
  Let us first assume that $A$ is an $\mathcal X$-filtered module and
  let $\{A_\alpha:\alpha \leq \sigma\}$ be an $\mathcal{X}$-filtration
  of $A$.  We are going to construct, by transfinite induction on
  $\alpha$, commutative diagrams
  \begin{displaymath}
    \commdiag{0 & \mapright & K_\gamma & \mapright^{g_\gamma} & \bigoplus_{\delta \leq
        \gamma}R^{(I_\delta)} & \mapright^{f_\gamma} & A_\gamma & \mapright & 0\cr
      & & \downarrow\lft{k_{\gamma\alpha}} & & \downarrow\lft{j_{\gamma\alpha}} & &
      \downarrow\rt{i_{\gamma\alpha}} & & \cr
      0 & \mapright & K_{\alpha} & \mapright_{g_\alpha} & \bigoplus_{\delta \leq
        \alpha}R^{(I_\delta)} & \mapright_{f_\alpha} & A_{\alpha} & \mapright & 0}
  \end{displaymath}
  in which $I_\gamma$ is a set, $j_{\gamma\alpha}$ and
  $i_{\gamma\alpha}$ are the inclusions and $k_{\gamma\alpha}$ is a
  monomorphism for each $\gamma < \alpha < \sigma$, satisfying the
  following two conditions:
  \begin{enumerate}
  \item If $\alpha<\sigma$ is a limit ordinal, then the corresponding row is
    the directed colimit of the previous ones.
  \item If $\alpha = \gamma +1$ is a successor ordinal, then
    $\textrm{Coker }k_{\gamma,\gamma+1} \in \mathcal{X}$.
  \end{enumerate}
  Note that the sequence constructed in step $\sigma$ will produce the
  desired presentation, since $K_\sigma$ is an $\mathcal{X}$-filtered
  module by construction.

  The case $\alpha = 0$ is trivial. Suppose now that $\alpha =
  \gamma+1\leq\sigma$ is a successor ordinal and that we have constructed the
  corresponding commutative diagrams for $\gamma$. Choose a presentation
  $\overline{f}_{\alpha}:R^{(I_\alpha)} \rightarrow
  \frac{A_{\alpha}}{A_\gamma}$ with $\ker \overline{f}_\alpha \in
  \mathcal X$ and, using the projectivity of $R^{(I_\alpha)}$, construct
  $h_\alpha:R^{(I_\alpha)}\rightarrow A_\alpha$ with
  $h_\alpha\pi_\alpha=\overline{f}_\alpha$, where $\pi_\alpha$ is the
  canonical projection. Let $f_{\alpha}:\bigoplus_{\delta \leq
    \alpha}R^{(I_\delta)} \rightarrow A_\alpha$ be the map induced by
  $f_\gamma$ and $h_\alpha$, $K_\alpha = \ker f_\alpha$. Call
  $\overline{K}_\alpha = \ker \overline{f}_\alpha$ and let
  $g_\alpha:K_\alpha \rightarrow \bigoplus_{\delta \leq
    \alpha}R^{(I_\delta)}$ and
  $\overline{g}_\alpha:\overline{K}_\alpha \rightarrow R^{(I_\alpha)}$ be
  the inclusions. Then using the Ker-Coker Lemma, we get the following commutative
  diagram with exact rows
  \begin{displaymath}
    \commdiag{ & & 0 & & 0 & & 0 & &\cr
      & & \downarrow & & \downarrow & & \downarrow & & \cr
      0 & \mapright & K_\gamma & \mapright^{g_\gamma} & \bigoplus_{\delta \leq
        \gamma}R^{(I_\delta)} & \mapright^{f_{\gamma}} & A_\gamma & \mapright & 0\cr
      & & \downarrow\lft{k_{\gamma\alpha}} & & \downarrow\lft{j_{\gamma\alpha}} & &
      \downarrow\rt{i_{\gamma\alpha}} & & \cr
      0 & \mapright & K_{\alpha} & \mapright^{g_\alpha} & \bigoplus_{\delta \leq
        \alpha}R^{(I_\delta)} & \mapright^{f_\alpha} & A_{\alpha} & \mapright &
      0\cr
      & & \downarrow & & \downarrow & & \downarrow & & \cr
      0 & \mapright & \overline{K}_\alpha &
      \mapright_{\overline{g}_\alpha} & R^{(I_\alpha)} &
      \mapright_{\overline{f}_\alpha} & \frac{A_\alpha}{A_\gamma} &
      \mapright & 0\cr
      & & \downarrow & & \downarrow & & \downarrow & & \cr
      & & 0 & & 0 & & 0 & & }.
  \end{displaymath}
  Set $k_{\delta \alpha} = k_{\delta \gamma}k_{\gamma \alpha}$,
  for each $\delta < \alpha$, and note that $\textrm{Coker
  }k_{\gamma\alpha} \in \mathcal{X}$, since $\ker \overline{f}_\alpha
  \in \mathcal{X}$.

  Finally, if $\alpha\leq\sigma$ is a limit ordinal, the directed
  colimit of the exact sequences for $\gamma < \alpha$ gives the
  desired sequence.

  Therefore, the result is true for $\mathcal X$-filtered modules. Let us now
  assume that $B$ is a direct summand of an
  $\mathcal{X}$-filtered module $A$. Then we have a commutative
  diagram with exact rows
  \begin{displaymath}
    \commdiag{0 & \mapright & K & \mapright & R^{(I)} & \mapright & A
      & \mapright & 0\cr
      & & \downarrow\lft{k} & & \downarrow\lft{1} & & \downarrow\rt{p} & & \cr
      0 & \mapright & K' & \mapright & R^{(I)} & \mapright & B &
      \mapright & 0}
  \end{displaymath}
  in which $K$ is $\mathcal{X}$-filtered and $p$ is a splitting
  epimorphism. Therefore, $k$ is also a splitting epimorphism
  and thus $K'$ is a direct summand of an $\mathcal{X}$-filtered module. In
  particular, it is also $\mathcal X$-filtered.\fin
\end{demostracion}

We can now prove the following result.

\begin{teorema}
  Let $\mathcal{S}$ be a set of modules containing the regular module $R$ and such that
  every module $S \in \mathcal S$ has a free presentation
  $\varphi:R^{(I)} \rightarrow S$ with $\ker \varphi \in \mathcal S$.
  Then the cotorsion pair cogenerated by $\mathcal{S}$ is hereditary.
\end{teorema}

\begin{demostracion}
  Follows from \cite[Theorem 2.2]{t} and propositions \ref{p:filtered}
  and \ref{p:hereditary}.\fin
\end{demostracion}

\begin{corolario} \label{p:Ahereditary} The cotorsion pairs $(\mathcal
  A, \mathcal B)$, $(\mathcal A^{<\infty},\mathcal B^{<\infty})$,
  $(\mathcal A_n, \mathcal B_n)$ and $(\mathcal A_n^{<\infty},\mathcal
  B_n^{<\infty})$, for each $n\in \mathbb{N}$, are hereditary.
\end{corolario}

\begin{demostracion}
  Let $n\in \mathbb{N}$. By Theorem \ref{t:cogenerated}, the cotorsion pair
  $(\mathcal A_n,\mathcal B_n)$ is cogenerated by the set $\{S \in
  \mathcal P_n: |S|<\kappa\}$, for some infinite regular cardinal
  $\kappa$. This set satisfies the hypothesis of the above
  theorem. Analogously, for each $S \in \mathcal P_n^{<\infty}$,
  we have that $|S| \leq |R|$. So the cotorsion pair $(\mathcal
  A_n^{<\infty},\mathcal B_n^{<\infty})$ is cogenerated by the set
  $\{S \in \mathcal P_n^{<\infty}:|S|<|R|^+\}$ and this set also
  satisfies the hypothesis of the previous theorem.\fin
\end{demostracion}

We get now the desired model structures on $\ch(R)$.

\begin{corolario}
  \begin{enumerate}
  \item The cotorsion pair $(\mathcal{A},\mathcal{B})$ (resp.
    $(\mathcal{A}_n,\mathcal{B}_n)$, for each $n\in \mathbb{N}$) induces a
    model category structure on $\ch(R)$ in which the weak
    equivalences are the homology isomorphisms, the cofibrations are
    the monomorphisms whose cokernels are in
    $\textrm{dg\,}\tilde{\mathcal{A}}$ (resp.
    $\textrm{dg\,}\tilde{\mathcal{A}}_n$), and the fibrations are the
    epimorphisms whose kernels are in $\textrm{dg\,}\tilde{\mathcal{B}}$
    (resp. $\textrm{dg\,}\tilde{\mathcal{B}}_n$).

  \item The cotorsion pair $(\mathcal{A}^{<
      \infty},\mathcal{B}^{<\infty})$ (resp.
    $(\mathcal{A}^{<\infty}_n,\mathcal{B}^{<\infty}_n)$, for each $n \in \mathbb{N}$) induces a model category structure on $\ch(R)$ in which
    the weak equivalences are the homology isomorphisms, the
    cofibrations are the monomorphisms whose cokernels are in
    $\textrm{dg\,}\tilde{\mathcal{A}}^{<\infty}$ (resp.
    $\textrm{dg\,}\tilde{\mathcal{A}}_n^{<\infty}$), and the fibrations
    are the epimorphisms whose kernels are in
    $\textrm{dg\,}\tilde{\mathcal{B}}^{<\infty}$ (resp.
    $\textrm{dg\,}\tilde{\mathcal{B}}_n^{<\infty}$).
  \end{enumerate}
\end{corolario}

\begin{demostracion}
 Apply Theorem \ref{t:ModelStructure}.\fin
\end{demostracion}

\section{Applications to the finitistic dimension}
\label{sec:appl-finit-dimens}

We finish this paper by applying our results to the calculus of
the finitistic dimensions of rings and algebras. From now on, we will assume that $\mathcal
F$ is the class of all projective modules and consequently, $\mathcal
C = R\Mod$. Therefore, $\mathcal P$ (resp. $\mathcal P^{<\infty}$) is
the class of all modules (resp. all modules in $R\modu$) with finite
projective dimension. Recall that the left big finitistic dimension of
$R$ is defined as
\begin{displaymath}
  \Findim(R) = \sup\{pd(P):P \in \mathcal P\}
\end{displaymath}
and the left little finitistic dimension of $R$ is
\begin{displaymath}
  \findim(R) = \sup\{pd(P):P \in \mathcal{P}^{<\infty}\}.
\end{displaymath}

Our goal is to characterize in Theorem \ref{t:Findim} when the big
finitistic dimension is finite. In order to do it, we will use the
following result which is an immediate consequence of the arguments
given in \cite[Lemma 2.2 and Proposition 2.3]{eego}.

\begin{proposicion}
  Let $\kappa$ be an infinite regular cardinal and $M$, a module such
  that both $M$ and each left ideal of $R$ are
  $<\kappa$-presented. Then $\displaystyle \ext^n\left(M,
    \lim_{\substack{\mapright\\\alpha < \kappa}}M_\alpha\right)
  \cong \lim_{\substack{\mapright\\\alpha <
      \kappa}}\ext^n(M,M_\alpha)$, for every $n\in \mathbb{N}$ and every
  directed system of morphisms $(M_{\alpha},f_{\alpha\beta})_{\alpha <
    \beta < \kappa}$.
\end{proposicion}

\begin{teorema}\label{t:Findim}
  Let $n$ be a natural number. The following assertions are
  equivalent:
  \begin{enumerate}
  \item $\Findim(R) \leq n$.

  \item $\pd(A) \leq n$ for every $A \in \mathcal A$.

  \item $\fibdim_{\mathcal B}(M) \leq n$ for every $M \in R\Mod$.

  \item $\fibdim_{\mathcal B}\left(R^{(R)}\right) \leq n$.
  \end{enumerate}
\end{teorema}

\begin{demostracion}
  1) $\Rightarrow$ 2) This follows from the facts that modules in
  $\mathcal{A}$ are direct summands of $\mathcal{P}$-filtered modules
  (see \cite[Theorem 2.2]{t}), the class $\mathcal P = \mathcal P_n$
  by hypothesis, and $\mathcal{P}_n$ is closed under direct summands
  and $\mathcal{P}_n$-filtrations (see \cite[Proposition 3]{a}).

  2) $\Leftrightarrow$ 3) Let $M$ be any module and take a fibrant replacement of $S(M)$ in $\ch(R)_{\mathcal M_\mathcal{A}^\mathcal{B}}$, $
    0 \to M \to B_0 \stackrel{{d_0}}{\mapright} B_1\stackrel{{d_1}}{\mapright} \cdots.$
  By Proposition \ref{p:CalculatingExt} we get that $\ext^1(A,B_{-n+1}B_\bullet)=\ext^{n+1}(A,M)=0$, , for every $A \in
  \mathcal A$ (since $\pd(A) \leq n$). This means that $B_{-n+1}B_\bullet \in
  \mathcal B$ and the result follows from the dual version of Proposition
  \ref{p:CalculatingDimensions}.

  3) $\Rightarrow$ 4) This is clear.

  4) $\Rightarrow$ 1) By Proposition \ref{p:CalculatingExt},
  $\ext^{m}\left(P,R^{(R)}\right)=0$ for every $P \in \mathcal{P}$ and
  every $m > n$. We deduce that $\ext^m\left(P,R^{(I)}\right) = 0$ for
  each $P \in \mathcal P$, each set $I$ and each $m > n$. Let us fix an $m > n$
  and a $P \in \mathcal P$ and let us denote by $\kappa = |R|^+$. Note
  that, by the proof of Theorem
  \ref{t:cogenerated}, we can assume that $P$ is $<\kappa$-presented.

  We are going to induct on the cardinality of $I$.  If $|I| \leq |R|$,
  the result is obvious. So assume that $|I| = \lambda > |R|$ and that
  $\ext^m\left(P,R^{(J)}\right) = 0$ for every set $J$ of cardinality strictly
  smaller than $\lambda$. Set $\mu_0 = \max \{\kappa,
  \textrm{cf}(\lambda)\}$ and note that $\mu_0$ and $P$ satisfy the
  hyphoteses of the above proposition. Moreover, since $\mu_0 \geq
  \textrm{cf}(\lambda)$, we can find an increasing sequence of
  ordinals, $\{\alpha_\nu:\nu < \mu_0\} \subseteq \lambda$, which is
  cofinal in $\lambda$. Therefore,
  \begin{displaymath}
    \ext^n\left(P,
      R^{(\lambda)}\right)
    \cong \lim_{\substack{\mapright\\\alpha <
        \mu_0}}\ext^n\left(P,R^{(\alpha)}\right) = 0
  \end{displaymath}
  by the induction hyphotesis.

  Finally, let us choose a module $P\in\mathcal P$ and let us suppose
  that $\textrm{pd}(P) = m>n$.  Let $M$ be any module and take a free
  presentation $0 \to K \to R^{(I)} \to M \to 0$
  of $M$. Then, applying the functor $\ext^{m}(P,\_)$ to this sequence, we get the exact
  sequence $\ext^m\left(P,R^{(I)}\right) \to \ext^m\left(P,M\right) \to \ext^{m+1}\left(P,K\right)$
  in which $\ext^m\left(P,R^{(I)}\right)=0$ by the previous considerations, and
  $\ext^{m+1}(P,K) = 0$ since $\textrm{pd}(P) = m$. But, as $M$ is
  arbitrary, this means that $\textrm{pd}(P) < m$. A contradiction
  that shows that $\textrm{pd}(P) \leq n$.\fin
\end{demostracion}

The above theorem can be easily improved when the ring $R$ is left perfect
and right coherent. Recall that, in this case, any direct product of
projective modules is projective (\cite[Theorem 3.3]{c}).

\begin{corolario}
  Let $R$ be a left perfect and right coherent ring, and $n$, a
  natural number. Then $\Findim(R) \leq n$ if and only if
  $\cofdim_\mathcal{B}(R) \leq n$.
\end{corolario}

\begin{demostracion}
  In this case, $R^{(I)}$ is a direct summand of $R^{I}$ for every
  index set $I$. So the hyphoteses of the corollary imply that
  $\ext^m\left(P,R^{(I)}\right) = 0$ for every $P \in \mathcal P$,
  every set $I$ and every $m > n$. Therefore, we can mimic the
  arguments used in the proof of 4) $\Rightarrow$ 1) in the above
  theorem in order to prove the result.\fin
\end{demostracion}

Out next result gives analogous results for the little finitistic
dimension.

\begin{teorema}
  Let $n$ be a natural number. The following assertions are
  equivalent:
  \begin{enumerate}
  \item $\findim(R) \leq n$.

  \item $\pd(A) \leq n$ for every $A \in \mathcal A^{<\infty}$.

  \item $\fibdim_{\mathcal B^{<\infty}}(M) \leq n$ for every $M \in
    R\Mod$.

  \item $\fibdim_{\mathcal B^{<\infty}}\left(R^{(R)}\right) \leq n$.
  \end{enumerate}
\end{teorema}

\begin{demostracion}
  The proof is analogous to that of Theorem \ref{t:Findim}. The only
  difference is that, in order to prove 4) $\Rightarrow$ 1), we need
  to use that $(\mathcal A^{<\infty},\mathcal B^{<\infty})$ is
  cogenerated by $\{A \in \mathcal P^{<\infty}:|A|<|R|^+\}$.\fin
\end{demostracion}

This Theorem can be improved for left coherent rings.

\begin{corolario}
  Let $R$ be a left coherent ring and $n$, a natural number. Then
  $\findim(R) \leq n$ if and only if
  $\fibdim_\mathcal{B^{<\infty}}(R) \leq n$.
\end{corolario}

\begin{demostracion}
  Since each module $P \in \mathcal P^{<\infty}$ is finitely presented
  and $R$ is left coherent, the functor $\ext^n(P,\_)$ commutes with
  direct limits for each $n\in \mathbb{N}$. Now, the fact that $\ext^m(P,R)
  = 0$ for each $m > n$ implies that $\ext^m\left(P,R^{(I)}\right) =
  0$ for each index set $I$. Finally, the proof of 4) $\Rightarrow$ 1)
  in Theorem \ref{t:Findim} gives the result.\fin
\end{demostracion}

\begin{corolario}\label{c:artinian}
  Let $R$ be a two-sided artinian ring.
  \begin{enumerate}
  \item The small finitistic dimension of $R$ is finite if and only if
    $\fibdim_{\mathcal B^{<\infty}}(R)$ is finite. Moreover, in this
    case, both dimensions do coincide.

  \item The big finitistic dimension of $R$ is finite if and only if
    $\fibdim_{\mathcal B}(R)$ is finite. Moreover, in this case,
    both dimensions do coincide.

  \item $\findim(R) = \Findim(R)$ if and only if
    $\fibdim_{\mathcal{B}^<\infty}(R) = \fibdim_\mathcal{B}(R)$.
  \end{enumerate}
\end{corolario}

\end{document}